\documentclass[11pt]{article}

\usepackage[cm]{fullpage}

\usepackage{amscd}
\usepackage{amsmath}
\usepackage{amsfonts}

\newcommand{\qed}{\hfill \ensuremath{\Box}}

\newcommand{\PP}{P}
\newcommand{\QQ}{\mathcal{Q}}
\newcommand{\FF}{\mathcal{F}}

\newcommand{\al}{\alpha}

\newcommand{\kahler}{K\"ahler }
\newcommand{\szego}{Szeg\"o }

\newcommand{\ecal}{\mathcal{E}}

\newcommand{\hcal}{\mathcal{H}}

\newcommand{\rcal}{\mathcal{R}}
\newcommand{\ocal}{\mathcal{O}}

\newcommand{\scal}{\mathcal{S}}

\def\us{{\underline S}}

\def    \t  {{\mathfrak t}}

\def    \Z  {{\mathbb Z}}
\def    \R  {{\mathbb R}}
\def    \C  {{\mathbb C}}

\newcommand{\CE}{{\mathcal E}}
\newcommand{\CR}{{\mathcal R}}
\newcommand{\CP}{\mathcal P}
\newcommand{\CPP}{\mathbf{CP}}

\newcommand{\minusone}{\sqrt{-1}}

\newcommand{\ddbar}{\sqrt{-1} \partial \overline{\partial}}

\newcommand{\ddt}[1]{\frac{\partial #1}{\partial t}}

\newcommand{\ddtt}[1]{\frac{\partial^2 #1}{\partial t^2}}

\numberwithin{equation}{section}

\begin{document}
\newcounter{remark}
\newcounter{theor}
\setcounter{remark}{0} \setcounter{theor}{1}
\newtheorem{claim}{Claim}
\newtheorem{theorem}{Theorem}[section]
\newtheorem{proposition}{Proposition}[section]
\newtheorem{lemma}{Lemma}[section]
\newtheorem{defn}{Definition}[theor]
\newtheorem{corollary}{Corollary}[section]
\newenvironment{proof}[1][Proof]{\begin{trivlist}
\item[\hskip \labelsep {\bfseries #1}]}{\end{trivlist}}
\newenvironment{remark}[1][Remark]{\addtocounter{remark}{1} \begin{trivlist}
\item[\hskip \labelsep {\bfseries #1
\thesection.\theremark}]}{\end{trivlist}}

\centerline{\Large \bf CONVERGENCE OF BERGMAN GEODESICS ON
$\mathbf{CP}^1$ \footnote{Research supported in party by National Science Foundation grants DMS-0604805  and DMS-0603850.} }

\bigskip
\bigskip

\centerline{\Large
\begin{tabular}{cc}
Jian Song & Steve Zelditch
\end{tabular} }

\bigskip

\centerline{\large
\begin{tabular}{c}
Department of Mathematics \\
Johns Hopkins University \\
\end{tabular}}

\bigskip

\bigskip
\bigskip

\section{Introduction}

This article is concerned with  geodesics in spaces of Hermitian
metrics of positive curvature on an ample line bundle $L \to X$
over a \kahler manifold.   Stimulated by a recent article of
Phong-Sturm \cite{PS}, we study the convergence as $N \to \infty$
of geodesics on the finite dimensional symmetric spaces $\hcal_N$
of Bergman metrics of `height $N$' to Monge-Amp\'ere geodesics on
the full infinite dimensional symmetric space $\hcal$ of
$C^{\infty}$ metrics of positive curvature. Our main result is
$C^{2}$ convergence of Bergman geodesics to Monge-Amp\'ere
geodesics in the  case of toric (i.e. $S^1$-invariant) metrics on
$\mathbf{CP}^1$. Although such metrics constitute the simplest
case of toric \kahler metrics,  the $\mathbf{CP}^1$ case already
exhibits much of the complexity of general toric varieties for the
approximation problem studied here. The general case will be
studied in \cite{SoZ}.

The convergence problem raised by Phong-Sturm \cite{PS} and
Arezzo-Tian \cite{AT}  belongs to the intensively studied program
initiated by Yau \cite{Y2} of relating the algebro-geometric issue
of stability to the analytic issue existence of canonical metrics
on holomorphic line bundles. In this program, metrics in $\hcal_N$
have a simple description in terms of algebraic geometry, while
metrics in $\hcal$ are `transcendental'. The approximation of
transcendental objects in $\hcal$ by `rational' objects in
$\hcal_N$ lies  at the heart of this program.

The reasons for studying Monge-Amp\'ere geodesics were laid out by
Donaldson in \cite{D1} (see also Mabuchi \cite{M} and  Semmes
\cite{S2}).   The  existence, uniqueness and regularity of such
geodesics is connected to existence and uniqueness  of metrics of
constant scalar curvature. Donaldson conjectured \cite{D2} that
there exist smooth Monge-Amp\'ere geodesics between any pair of
metrics $h_0, h_1 \in \hcal$. The best general result, due to Chen
\cite{Ch} and Chen-Tian \cite{CT}, is the existence of a unique
$C^{1,1}$ geodesic $h_t$ joining $h_0$ to $h_1$. This result is
sufficient to prove uniqueness of extremal metrics. In the case of
toric varieties, the much stronger result is  known that the
geodesic between any two metrics is $C^{\infty}$ \cite{G}. In
fact, the Monge-Amp\'ere equation can be linearized by the
Legendre transform and the symmetric space is flat.

But in general, solutions of  the  Monge-Amp\'ere equation are
difficult to analyze.  It seems rather remarkable that one can
study solutions of the homogeneous Monge-Amp\`ere equation by
means of `algebro-geometric approximations'. It has been
proved by Phong-Sturm \cite{PS} (see also \cite{B})  that Bergman
geodesics, which are orbits of one-parameter subgroups of $GL(d_N
+ 1, \C)$ between two Bergman metrics, converge uniformly to a
given Monge-Amp\'ere geodesic for a general ample line bundle over
a \kahler manifold.

The question we take up in this article and in \cite{SoZ} is
whether Bergman geodesics converge to Monge-Amp\`ere geodesics in
a stronger sense. Convergence in $C^2$ is especially interesting
since it implies that the curvatures and moments maps for the
metrics along the Bergman geodesic converge to those along the
Monge-Amp\`ere geodesics.  In this article and in the subsequent
article \cite{SoZ}, we study this problem
 for toric hermitian line bundles over toric \kahler
manifolds. In this setting,  the \kahler potentials $\varphi_N(t,
z)$ of the Bergman metrics along the geodesic have  relatively
explicit formulae (see \ref{intphi2})  resembling the free energy
of a discrete quantum statistical mechanical model. Convergence in
$C^0$ of the \kahler potential as $k \to \infty$ is analogous to
uniform convergence of the free energy in the thermodynamic limit,
while convergence of derivatives is related to absence of phase
transitions (cf. \cite{E}, II.6).

To state our results, we will need some notation. Let  $L \to X$
be an ample holomorphic  line bundle and denote by  $H^0(X, L^N)$
the holomorphic sections of the $N^{th }$ power $L^N \to X$ of $L$.
 Given a basis ${\mathcal S_N} = \{ S_0, \dots,
S_{d_N}\}$ we define the associated holomorphic embedding
\begin{equation}\label{KODAIRA} \Phi_{{\mathcal S}_N}: X \to \mathbf{CP}^{d_N},\;\;\; \Phi_{{\mathcal S}}(z) = [S_0(z), \dots, S_{d_N}(z) ].
\end{equation}
 We define the space of Bergman metrics by $$ \mathcal{H}_N = \{ \frac{1}{N} \Phi_{\mathcal{S}}^* \omega_{FS} ~|~  \mathcal{S}  {\rm ~is ~ a ~ basis ~ of ~   }
 H^0(X, L^N) \}, $$ where $\omega_{FS}$
is the Fubini-Study metric on $\mathbf{CP}^{d_N}$. Since $U(d_N+1)$
is the isometry group of $\omega_{FS}$, ${\mathcal H}_N$ is the
symmetric space $ GL(d_N+1, \C)/U(d_N+1,\C)$.

Metrics in  $\hcal_N$ are defined by an essentially
algebro-geometric construction and are somewhat analogous to
rational numbers.  A basic fact is that the union
$$ \bigcup_{N =1}^{\infty} \hcal_N \; \subset \; \hcal$$
of Bergman metrics  is dense in the $C^{\infty}$ topology  in the
space $\hcal$ of all $C^{\infty}$ K\"ahler metrics (see \cite{T,
Z}) of positive curvature. Indeed,   for each $N$ we have a map
\begin{equation} \label{BCALK} \scal_N: \hcal\to \hcal_N, \;\; h
\to \; h_N = (\Phi_{\mathcal{S}_N}^* h_{FS} )^{1/N},\;\; \mathcal{S}_N = \mbox{an
orthonormal basis for}\;\; h.
\end{equation}  The metric $h_N$ is independent of the
choice of orthonormal basis, and $h_N\to h$ in $C^{\infty}$.

Now let us compare Monge-Amp\'ere geodesics and Bergman geodesics.
We let $h_0, h_1$ be any two  hermitian metrics on $L$ in the
class $\hcal$ and write  $h_{\varphi} = e^{-\varphi} h$ relative to a
fixed metric $h$ with curvature form $\omega=Ric(h)$. Thus, we have
an isomorphism \begin{equation} \label{HCALDEF}  \hcal\ = \
\{\varphi\in C^{\infty} (X) ~|~ \omega_\varphi\ = \
\omega+\ddbar \varphi>0\ \} \ .
\end{equation}
  We may then
identify the  tangent space $T_\varphi\hcal$ at $\varphi\in\hcal$  with
$C^{\infty}(X)$. We define a Riemannian metric on $\hcal$ as
follows: let  $\varphi\in \hcal$ and let $\psi\in T_\varphi\hcal \simeq
C^{\infty} (M)$ and define
\begin{equation} \label{metric} ||\psi||^2_{\varphi}\ = \ \int_M |\psi|^2\
\omega_{\varphi}^n\ \ . \end{equation} With this Riemannian metric,
$\hcal$ is an infinite dimensional negatively curved symmetric
space. By \cite{M, S1, S2, D1}, the geodesics of $\hcal$ in this
metric are the paths $\varphi_t$ which satisfy the equation
\begin{equation} \ddot\varphi-|\partial\dot\varphi |_{\omega_{\varphi}}^2=0.
\end{equation}
This may be interpreted as a Monge-Amp\'ere equation \cite{S1,
D1}.

Geodesics in $\hcal_N$ with respect to the symmetric space metric
are given by one-parameter subgroups $e^{t A}$ of $GL(d_N+1, \C)$.
That is, let  $\sigma\in GL(d_N+1,\C)$ be the change of basis matrix
defined by $\sigma\cdot\hat\us^{(0)}=\hat\us^{(1)}$. Without loss
of generality, we may assume that $\sigma$ is diagonal with
entries $e^{\lambda_0},...,e^{\lambda_{d_N}}$ for some
$\lambda_j\in\R$. Let $\hat\us^{(t)}=\sigma^t\cdot\hat\us^{(0)}$
where $\sigma^t$ is diagonal with entries $e^{\lambda_jt}$.
Define
\begin{equation}  h_N(t,z)=h_{\hat\us^{(t)}}(z)=h(z)e^{-\varphi_N(t,z)}. \end{equation}  Then
$h_N(t,\cdot)$ is the smooth geodesic in $GL(d_N+1,\C)/U(d_N+1,\C)$ joining
$h_N(0,\cdot)$ to $h_N(1, \cdot)$. Explicitly,  we have
\begin{equation} \label{intphi}
   \varphi_N(t,z)\ = \ {1\over N}
\log\left(\sum_{j=0}^{d_N} e^{2\lambda_jt}|\hat S_j^{(0)}|^2_{h^N}(z)
\right).
\end{equation}

Thus, the problem is the convergence of $h_N(t,\cdot) \to h(t,\cdot)$ or
equivalently of $\varphi_N(t,\cdot) \to \varphi(t,\cdot)$. The following general result is proved in \cite{PS}.

\begin{theorem} The Bergman geodesics uniformly converge to the Monge-Amp\`ere geodesic in the sense that
\begin{equation}
\varphi_t(z)=\lim_{ k \to \infty }[\sup_{N\geq k} \varphi_N(t, \cdot)]^*(z),
\end{equation}
where, for any bounded function $f: [0,1]\times X\to \mathbf{R}$, the upper envelope of $f$ is defined by $f^*(x_0)=\lim_{\epsilon\to0} \sup_{|x-x_0|<\epsilon} f(x)$.
\end{theorem}
As mentioned above, our  goal here and in \cite{SoZ} is to study
the degree of convergence of these geodesics in the case of  toric
hermitian metrics on a toric line bundle $L \to X$. We define the
space
 $\hcal_{T}$ to be the subspace of $\hcal$ of hermitian metrics for which
 $\varphi$ is invariant under the underlying real torus $T = (S^1)^n$ action.

In the case of $\mathbf{CP}^1$, we may assume $L = \ocal(1)$ and
an orthogonal basis of holomorphic sections of $L^N = \ocal(N)$ is
given in an affine chart by the monomials $z^{\alpha}, \alpha = 0,
\dots, N$. A toric hermitian metric  is entirely encoded in the
set of  $L^2$ squared norms $\mathcal{Q}_h^N(\alpha) =
||z^\alpha||_{h^N}^2$ of the monomials with respect to powers
$h^N$ of the Hermitian metric $h$ (cf. Definition
\ref{NORMINGCONSTANTS}). Then (\ref{intphi}) takes the form
\begin{equation} \label{intphi2}
   \varphi_N(t,z)\ = \ {1\over N}
\log\left(\sum_{\alpha \in \frac{1}{N} \Z \cap [0, 1]} \frac{|z^{N
\alpha}|_{h_0}^2}{(\mathcal{Q}_{h_0}^N(\alpha))^{1-t}
(\mathcal{Q}_{h_1}^N(\alpha))^t}  \right).
\end{equation}
If we write $|z|^2 = e^{\rho}$, we see the resemblance to the free
energy of a quantum statistical model with states parameterized by
lattice points in $[0, N]$ \cite{E} (\S 7). The main result of
this article is:

\begin{theorem} \label{main} On  $\mathbf{CP}^1$, the Bergman geodesics converge to the toric Monge-Amp\`ere geodesic uniformly
\begin{equation}
\lim_{N \to \infty} \varphi_N(t, z) = \varphi_t (z),
\end{equation}
uniformly in the $C^{2}$ topology on $[0,1]\times \CPP^1$.
\end{theorem}

 As will be seen in the proof, most of the complications
concern the joint asymptotics in the $(N, \alpha)$ parameters of
the {\it norming constants} $\mathcal{Q}_h^N(\alpha)$  near the
boundary of the `moment polytope' $[0, 1]$. The essential
simplification in $\mathbf{CP}^1$ over higher dimensional toric
varieties is that the approach to the boundary is much simpler for
an interval than for the possible convex Delzant polytopes in
higher dimensions. Otherwise, the case of $\mathbf{CP}^1$ already
exhibits much of the complexity of the general case.  In
\cite{SoZ}, we study the $C^2$ convergence problem in all
dimensions.

 Our analysis of the norming constants
builds on the work of  \cite{STZ1}, and may have an independent
interest, since the norming constants  determine a toric metric.
For instance, in \cite{D4} and elsewhere, numerical methods for
approximating extremal \kahler metrics on toric varieties are also
based on the study of norming constants.  It would be interesting
to generalize the results on norming constants to higher
dimensions. The subsequent article \cite{SoZ} involves quantities
which are in a sense dual to norming constants and does not
directly provide information on norming constants.

\section{Preliminaries}

Although we primarily study $\mathbf{CP}^1$ in this article, we
set the scene for toric varieties in arbitrary dimensions.  Let
$(X,\omega, \tau)$ be a compact toric manifold of complex
dimension $n$ and
$$\tau: T^n\rightarrow Diff(X, \omega)$$ an effective Hamilton
action of the standard real $n$-torus $T=(S^1)^n$. Let $\pi$ be
the moment map associated to the toric K\"ahler metric $\omega$
\begin{equation} \label{MM} \pi: X\rightarrow \mathbf{R}^n. \end{equation} The image $P$ of
$\pi$ is a Delzant polytope,  defined by a set of linear
inequalities given by
$$\langle x, v_r\rangle\geq \lambda_r, ~~~r=1, ..., d, $$
where $v_r$ is an inward-pointing normal to the $r$-th
$(n-1)$-dimensional face of $P$. Define the affine functions $l_r:
\mathbf{R}^n\rightarrow \mathbf{R}$ by
$$l_r(x)=\langle x, v_r\rangle-\lambda_r.$$

Fix a toric polarization $L$ on $X$ with $[L]=[\omega]$. Let
$$\mathcal{H}=\{ h~|~ h ~\textnormal{is
a smooth $T$-invariant hermitian metric on $L$ such that $Ric(h)>0$}\}.$$
Fix $h \in \mathcal{H}$ and let $\omega = Ric(h)$, then
$$\mathcal{H} \cong \{ \varphi\in C^{\infty} (X)~|~\textnormal{$\varphi$ is $T$-invariant and ~}\omega_{\varphi}
=\omega + \ddbar\varphi>0\}.$$
Hence the hermitian metric $h_{\varphi}  \in\mathcal{H}$ and the
plurisubharmonic potential $\varphi\in \mathcal{H}$ are related by
$$h_{\varphi}=h_0 e^{-\varphi}.$$
The $L^2$-metric on $\mathcal{H}$ is given by
$$||\psi||_{\omega_{\varphi}}^2=\int_X
|\psi|^2\omega_{\varphi}^n$$
for any $\psi\in C^{\infty}(X)$.

 For any $\varphi_0$ and $\varphi_1\in
\mathcal{H}$, the geodesic $\varphi_t$ joining $\varphi_0$ and
$\varphi_1$ in $\mathcal{H}$ is defined by

\begin{equation}
\frac{\partial^2\varphi_t}{\partial t^2}=\left| \partial
\dot\varphi_t\right|_{\omega_{\varphi_t}}^2.
\end{equation}

From a complex geometric viewpoint, the complex torus
$(\mathbf{C}^*)^n$ acts on $X$ with an open orbit, and $X$ may be
viewed as  a compactification
 of $(\mathbf{C}^*)^n$. On  the open orbit, we denote the standard  holomorphic coordinates
by $(z_1, ... , z_n)$. We also define the real coordinates $
\rho_j=\log |z_j|^2$, $j=1, ... , n$. Then a toric \kahler form
has a $T$-invariant \kahler potential $u$ on the orbit
defined by $\omega= \sum_{i, j=1, ..., n} \sqrt{-1} \frac{\partial ^2 u}{\partial
z_i
\partial \overline{z}_j} d z_i\wedge d \overline{z}_j>0$. Since $u$ is
$T$-invariant, it can be considered as a function in
$\rho=(\rho_1, ... , \rho_n)$ on $\mathbf{R}^n$ and it is convex
on $\mathbf{R}^n$. We then define $U(\rho)= u(z)$ on $\mathbf{R}^n$.

The Legendre transform of $U$ defines the symplectic potential $G$
of $\omega$, a convex function on $P^{\circ}$.  That is,

$$G(x)= \langle x, \rho \rangle -U(\rho)$$
with $x=\nabla U (\rho)\in P \subset \mathbf{R}^n$  given by the
moment map. It has the same  singularities at the boundary
$\partial P$ as the symplectic reference potential
\begin{equation}
G_{P}(x)=\sum_{r=1}^{d}l_r(x)\log l_r(x).
\end{equation}
$G_P$ induces a smooth hermitian metric $h_P$ on $L$ with
$Ric(h_P)=\ddbar u_P$ on $(\mathbf{C}^*)^n$\ and $u_P$ being the
Legendre transform of $G_P$. For background, we refer to \cite{A,
D4,Gu}.

The following theorem is proved by Guan \cite{G}.
\begin{theorem}
Let $h_t$ be the smooth geodesic joining $h_0$ and $h_1\in \mathcal{H}$ for $t\in [0,1]$. The
corresponding symplectic potential $G_t$ is given by
\begin{equation}
G_t(x)=G_{P}(x)+f_t(x)
\end{equation}
where $f_t$ is a smooth function on $\mathbf{R}^n$ with $\nabla^2
G_t>0$ on $P^{\circ}$. Furthermore,

\begin{equation}
f_t(x)=(1-t)f_0(x)+tf_1(x).
\end{equation}

\end{theorem}
Hence the geodesic of the symplectic potentials is linear. A very
simple proof (cf. \cite{SoZ}) is simply to push forward the energy
functional defining the Monge-Amp\`ere geodesics to the polytope
and observe that it becomes the Euclidean energy functional there.

\begin{defn} \label{NORMINGCONSTANTS}
For any lattice point $N\alpha\in NP\cap\mathbf{Z}^n$, we define the
$L^2$ norm of $S^N_{\alpha}\in H^0(X, L^N)$ with respect to $h_t$
by

\begin{equation}
\mathcal{Q}_t^N(\alpha)=\int_X |S^N_{\alpha}|^2_{h_t^N} \omega_t^n
\end{equation}
where $\omega_t=Ric(h_t)$ and $h_t^N$ the $N^{th}$-power of $h_t$. We also define
$\mathcal{Q}_{P}^N(\alpha)$ with respect to $h_P$ by

\begin{equation}
\mathcal{Q}_{P}^N(\alpha)=\int_X |S^N_{\alpha}|^2_{h_{P}^N}
 \omega_{P}^n
\end{equation}
where $\omega_P$ is the toric K\"ahler form given by the symplectic potential $G_P$. Since $S^N_{\alpha}$ is defined as the monomial $z^{\alpha}$ on $(\mathbf{C}^*)^n$,  $\QQ_t^{N}(\alpha)$ and $\QQ_P^N(\alpha)$ can also be defined for all $\alpha\in P$.

\end{defn}

Phong and Sturm \cite{PS} introduce the $GL(d_N+1, \mathbf{C})$ geodesics in
the space of Bergman metrics to approximate the Monge-Ampere
geodesic $\varphi_t$.

\begin{defn} \label{ECAL}
We define $\mathcal{E}_N(t, z)$ by
\begin{equation}
\mathcal{E}_N(t, z)=\sum_{N\alpha\in
NP\cap\mathbf{Z}^n}
\frac{|S^N_{\alpha}|^2_{h_t^N}}{(\mathcal{Q}_0^N(\alpha))^{1-t}
(\mathcal{Q}_1^N(\alpha))^t},
\end{equation}
and the Szeg\"o kernel $\Pi_N$ with respect to $h_t$ by

\begin{equation}
\Pi_N(t, z)=\sum_{N\alpha\in N P\cap\mathbf{Z}^n}
\frac{|S^N_{\alpha}|^2_{h_t^N}}{Q_t^N(\alpha)}.
\end{equation}
\end{defn}

\begin{defn} \label{QVSQQ}

We also define for $\alpha \in P$

\begin{enumerate}

\item the norming constants

\begin{equation} \label{Q}
Q_t^N(\al)=\QQ_t^N(\al)
e^{-NG_t(\al)},~~~Q_{\PP}^N(\al)=\QQ_{\PP}^N(\al)
e^{-NG_{\PP}(\al)},
\end{equation}

\item the norming constants \begin{equation}
q_t^N(\alpha)=\frac{Q_t^N(\alpha)}{Q_{\PP}^N(\alpha)},~~~
\mathcal{R}_t^{N}(\alpha)=\frac{q_t^N(\alpha)}{(q_0^N(\alpha))^{1-t}(q_1^N(\alpha))^t}
=\frac{Q_t^N(\alpha)}{(Q_0^N(\alpha))^{1-t}(Q_1^N(\alpha))^t},
\end{equation}

\item the norm squares of the normalized monomial sections

\begin{equation}\label{PCAL}
\mathcal{P}^N_{\alpha} (t, z)=\frac{|S^N_{\alpha}|^2_{h_t^N}(z)}{\mathcal{Q}_t^N(\alpha)}.
\end{equation}

\end{enumerate}

\end{defn}

\begin{lemma}\label{RCAL}

\begin{equation}
\mathcal{E}_N (t, z)  =\sum_{N\alpha\in
N \PP\cap \mathbf{Z}^n}\mathcal{R}_t^N(\alpha)\mathcal{P}^N_{\alpha}(t, z).
\end{equation}

\end{lemma}

\begin{proof} Straightforward calculation shows that
\begin{eqnarray*}
&&\frac{q_t^N(\alpha)}{ (q_0^N(\alpha))^{1-t}(q_1^N(\alpha))^{t}}\\
&=&
\frac{Q_t^N(\alpha)}{ (Q_0^N(\alpha))^{1-t}(Q_1^N(\alpha))^{t}}\\
&=& e^{N \left( (1-t)G_0(\al)+tG_1(\al)-G_t(
\al)\right)}\frac{\QQ_t^N(\alpha)}{ (\QQ_0^N(\alpha))^{1-t}(\QQ_1^N(\alpha))^{t}}\\
&=& \frac{\QQ_t^N(\alpha)}{
(\QQ_0^N(\alpha))^{1-t}(\QQ_1^N(\alpha))^{t}}.
\end{eqnarray*}
The last equality follows from the geodesic equation
$G_t(x)=(1-t)G_0(x)+tG_1(x)$.

\qed
\end{proof}

\section{\label{OINTQ} Joint $(N, \alpha)$ asymptotics of the  norming constants for metrics on $\mathbf{CP}^1$}

We first give a useful formula for the norming constants
$Q_t^N(\alpha)$ (\ref{Q}) which is valid on any toric variety, and
then we use it in the case of $\mathbf{CP}^1$ to determine joint
$(N, \alpha)$ asymptotics.

\begin{lemma}\label{phase}
The  norming constants $Q_t^N(\alpha)$ and $Q_P^N(\alpha)$ in
Definition \ref{QVSQQ} for $\alpha \in P$ are given on any toric
variety by

\begin{eqnarray}
&&Q_t^N(\alpha)=(2\pi)^n\int_{\PP} e^{-NF_{t,\alpha}(x)}dx\nonumber\\
&&Q_P^N(\alpha)=(2\pi)^n\int_{\PP} e^{-NF_{P,\alpha}(x)}dx,
\end{eqnarray}
where the phase functions $F_{t,\alpha}(x)$ and $F_{P,\alpha}(x)$
are defined by
\begin{equation}\label{PHASEDEF} \left\{\begin{array}{lll}
F_{P,\alpha}(x) & = & \langle x-\alpha, \nabla
G_{\PP}(x)\rangle-(G_{\PP}(x)-G_{\PP}(\alpha))\\
&&\\
F_{t,\alpha}(x) & = & \langle x-\alpha, \nabla
G_t(x)\rangle-(G_t(x)-G_t(\alpha)).
\end{array} \right.\end{equation}

\end{lemma}

\begin{proof}  Let $z=(z_1, ... , z_n)\in  \left( \mathbf{C}^* \right)^n$ and $\rho=(\rho_1, ..., \rho_n) \in \mathbf{R}^n$ with $\rho_j=\log |z_j|^2$ for $j =1, ... , n$. We suppose that the K\"ahler form for $g_t$ is given by $\sum_{i, j=1, ..., n } \sqrt{-1} \frac{\partial^2 u_t}{\partial z_i \partial \overline{z}_j} dz_i\wedge d\overline{z}_j$, where $u_t(z)$ is the K\"ahler potential for the toric K\"ahler metric $g_t$ on $ \left( \mathbf{C}^*\right)^n$. Let $U_t (\rho)=u_t (z)$ and $\pi_t=\nabla U_t: \mathbf{R}^n \rightarrow P$ be the moment map associated to $g_t$. Then the symplectic potential $G_t$ on $P$ for $g_t$ is given by the following Legendre transform

$$ G_t(x)= \langle x, \rho \rangle - U_t(\rho) $$
with $x=\nabla U_t(\rho)\in P \subset \mathbf{R}^n$. Also $U_t(\rho)$ can be recovered from $G_t(x)$ by the following inverse Legendre transform

$$ U_t( \rho )= \langle x, \rho \rangle - G_t(x)$$
with $\rho= \nabla G_t(x)$. Also $\pi_t^* ( dx_1... dx_n ) = \det \left(  \frac{\partial^2 U_t}{\partial \rho_i \partial \rho_j}  \right)   d\rho_1...d\rho_n$.

\begin{eqnarray*}
Q_t^N(\alpha) &=& (\minusone)^n  \int_{ \mathbf{C}^n}  |z|^{N\alpha} e^{-Nu_t(z)-N G_t(\alpha)} \det \left( \frac{ \partial ^2 u_t}{\partial z_i \partial \overline{z}_j}  \right) d z_1\wedge d \overline{z}_1 \wedge ... \wedge     d z_n\wedge d \overline{z}_n\\
&=& (2\pi)^n \int_{\mathbf{R}^n} e^{N(\langle \alpha, \rho \rangle -U_t(\rho))-N G_t(\alpha)} \det \left(  \frac{\partial^2 U_t}{\partial \rho_i \partial \rho_j}  \right)   d\rho_1...d\rho_n\\
&=& (2\pi)^n \int_P e^{N( \langle \alpha, \nabla G_t(x) \rangle - (\langle x , \nabla G_t(x) \rangle- G_t(x))- G_t(\alpha))} dx_1... d x_n \\
&=& (2\pi)^n \int_P e^{-N F_{\t, \alpha}(x)} dx.
\end{eqnarray*}

The same argument gives the integral formula for $Q_P^N(\alpha)$.

\qed
\end{proof}

We now specialize to the case of  $\mathbf{CP}^1$, where:

\begin{itemize}

\item  $P = [0, 1]$ and   $G_P(x)=x\log x+(1-x)\log(1-x)$;

\item  For $\alpha \in \frac{1}{N} \Z \cap P$, $\QQ_{\PP}^N(\al) =
{N \choose N \alpha}^{-1},$ and $Q_P^N(\alpha) = 2\pi {N \choose N
\alpha}^{-1} e^{- N \left(\alpha \log \alpha + (1 - \alpha) \log
(1 - \alpha)\right)}.$

\item  The geodesic of the symplectic potentials $G_t(x)$ is
\begin{equation*}
G_t(x)=G_P(x)+f_t(x)
\end{equation*}
where $f_t(x)=(1-t)f_0(x)+t f_1(x)$ is a smooth function on
$\mathbf{R}$
 such that
\begin{equation}\label{GTPOS}
\frac{d^2}{dx^2}G_t(x)>0.
\end{equation}

\end{itemize}

 In fact,  because $G''_t(x)$ has poles of order $1$ at $0$ and
$1$, we have:
\begin{lemma}\label{bound1} There exists a constant $\Lambda>0$ such that for any
$t\in [0,1]$ and $x\in (0,1)$

\begin{equation}
x(1-x)G''_t (x)>\Lambda,  ~~~~  x(1-x)G''_P(x)>\Lambda.
\end{equation}

\end{lemma}

We also evaluate:
\begin{equation}  \label{FTALPHA} \left\{ \begin{array}{lll}
F_{P,\alpha}(x)&=&-\alpha\log
x-(1-\alpha)\log(1-x)+\alpha\log\alpha+(1-\alpha)\log(1-\alpha),\\
&& \\
F_{t,\alpha}(x)
%
%
&=&G_P(\alpha)-\alpha\log x-(1-\alpha)\log(1-x)
+(x-\alpha)^2f_{t, \alpha}(x),\\
\end{array} \right.\end{equation}
where
$f_{t,\alpha}(x)=-\frac{f_t(x)-f_t(\alpha)-f'_t(x)(x-\alpha)}{(x-\alpha)^2}$
with $f_{t, \alpha}(\alpha)=\frac{1}{2} f_t''(\alpha)$. It is easy
to check by Taylor expansion that $f_{t,\alpha}(x)$ is smooth in
$x$ and $t$.

%


We now consider the  joint asymptotics in $(N, \alpha)$ of the
norming constants. Our main result, Theorem \ref{main1}, is a
comparison of  the joint asymptotics of a metric norming constant
(\ref{Q}) to  the canonical norming constants $Q_P^N(\alpha)$. The
joint asymptotics of the latter   can be derived from known
(elementary) results on binomial coefficients, and we begin by
recalling the relevant background.

The joint asymptotics of binomial coefficients ${N \choose m}$ in
$(N, m)$ and the closely related canonical norming constants
$Q_P^N(\alpha)$ have several regimes accordingly as $\alpha$
belongs to an `interior region' or a `boundary region'. First let
us consider the  `left interior,' where $\alpha \in
[\frac{1}{N^{3/4}}, 1 - \frac{1}{N^{3/4}}]$. The standard Sterling
asymptotics for factorial and  binomials applies in the region and
gives
\begin{equation} \label{BINOMIAL}
{N \choose N \alpha} \sim \;\;\; \frac{1}{\sqrt{2\pi N \alpha(1 -
\alpha)}} \; e^{- N (\alpha \log \alpha + (1 - \alpha) \log (1 -
\alpha)},
\end{equation}
and more precisely the asymptotics \begin{equation} \label{QPINT}
Q_P^N(\alpha) = 2\pi {N \choose N \alpha }^{-1}  e^{- N G_P
(\alpha)} = 2\pi \sqrt{(2 \pi) N \alpha (1 -  \alpha)}  \exp
\left( O\left(\frac{1}{N \alpha} + \frac{1}{N - N \alpha} \right)
\right).
\end{equation}  We observe that the asymptotics are highly non-uniform
as $\alpha \to 0$ or $\alpha \to 1$.

In the left `boundary region' $\alpha \in [0, \frac{1}{N^{3/4}}],$
we cannot use Stirling's formula up the boundary  and rather use
that
$${N  \choose  m} = A(N, m)
\frac{N^{m}}{m!},\;\; \mbox{with}\; A = \Pi_{j = 1}^{m - 1} (1 -
\frac{j}{N}). $$ Using that $\ln A = \sum_{j = 1}^{m - 1} \ln (1 -
\frac{j}{N}), $ and $\ln (1 - x) \sim -x$ one has
$$\sum_{j = 1}^{N\alpha - 1} \ln (1 - \frac{j}{k}) \sim \sum_{j = 1}^{N\alpha - 1}
- \frac{j}{N} \sim \frac{(N\alpha)^2}{2 N} = o(1) $$ if $N \alpha
= o(\sqrt{N}). $
 It follows that  if $(N \alpha) = o(\sqrt{N}) $,
then ${N \choose N \alpha} \sim \frac{N^{N \alpha}}{(N \alpha)!}$,
and further that \begin{equation} \label{QPBD} \begin{array}{lll}
2\pi \left(Q_{\PP}^N(\al)\right)^{-1} &
=& {N \choose N \alpha} (\frac{N \alpha}{N})^{N \alpha} (1 - \frac{N \alpha}{N})^{N - N \alpha} \\ && \\
& = & \frac{(N \alpha)^{N\alpha}}{(N \alpha)!} (1 - \frac{1}{N})
(1 - \frac{2}{N}) \cdots (1 - \frac{N \alpha}{N})(1 -
\frac{N \alpha}{N})^{N - N \alpha} \\ && \\
& = & (1 - \frac{(N \alpha)}{N})^{N} \frac{(N \alpha)^{(N
\alpha)}}{(N \alpha)!} (1 - \frac{1}{N}) (1 - \frac{2}{N}) \cdots
(1 - \frac{N \alpha}{N})(1 - \frac{N \alpha}{N})^{- N \alpha}.
\end{array} \end{equation}

We record the following:

\begin{lemma}\label{LBCQ}

There exists a constant $C>0$ such that for all $\alpha \in [0,1]\cap
\frac{1}{N}\mathbf{Z}$

\begin{equation}\label{polybound}
Q_P^N(\alpha)\geq C.
\end{equation}

\end{lemma}

\begin{proof} In the interior region, (\ref{QPINT}) implies the lower bound
$Q_P^N(\alpha) \geq C N^{1/8}$.  In the boundary region, we can
continue to  use Stirling's formula as long as $N \alpha \to
\infty$  to obtain
$${N \choose N \alpha} \sim \left(\frac{N e}{N \alpha} \right)^{N \alpha} (2 \pi N
\alpha)^{-1/2} \implies  {N \choose N \alpha} \left(\frac{N
\alpha}{N}\right)^{N \alpha} (1 - \frac{N \alpha}{N})^{N - N
\alpha} \sim (2 \pi N \alpha)^{-1/2},$$ so $Q_P^N(\alpha)  \to
\infty$ there as well. If $N \alpha \leq K$ then the exact formula
(\ref{QPBD}) gives  positive upper bound independent of $N$. We
note that it equals $1$ when $\alpha = 0$.

\end{proof}

 We now turn to general metrics. The following comparison
 inequality
 is the principal technical tool in the proof of
$C^2$ convergence of the geodesics (see Definition
\ref{QVSQQ}).

\begin{theorem}\label{main1}
There exists a constant  $C>0$ such that for all integer $N >0$, $\alpha\in P$ and $t \in [0,1]$
\begin{equation}
\frac{1}{C}\leq q_t^N(\alpha)\leq C.
\end{equation}
Furthermore, if we let $\pi_t$ and $\pi_{\PP}$ be the moment maps
associated to the toric K\"ahler metrics $g_t$ and $g_{\PP}$ and define
$$\Omega_t(\al)=\left(\frac{\det
\nabla^2G_{\PP}(\al)}{\det
\nabla^2G_t(\al)}\right)^{\frac{1}{2}},$$  then $\Omega_t(\alpha)$
extends to a continuous function on $\PP$ and
\begin{equation}
\lim_{N\rightarrow\infty} q_t^N(\al)=\Omega_t(\al)
\end{equation}
uniformly for  $\alpha \in \PP$.

\end{theorem}

Indeed,
$\Omega_t(\alpha)=\left(\frac{\det \nabla^2 \Phi_t(\pi_t^{-1}(\al))}{\det
\nabla^2 \Phi_P (\pi_{\PP}^{-1}(\al))}\right)^{\frac{1}{2}}$ is the ratio
of the volume forms of K\"ahler metrics $g_t$ and $g_P$ on $(\mathbf{C}^*)^n$,
 although $\pi_t^{-1}(\al)$ and $\pi_P^{-1}(\al)$ do not necessarily coincide.

 The following corollaries play an important role in the proof of
 the main result.

\begin{corollary} \label{maincor1}
There exists a constant  $C>0$ such that for all integer $N >0$,
$\alpha\in [0, 1]$ and $t \in [0,1]$, the ratios $\mathcal{R}$ of
Defintion \ref{QVSQQ} satisfy
\begin{equation}
\frac{1}{C}\leq \mathcal{R}_t^N(\alpha)\leq C.
\end{equation}
Furthermore,
\begin{equation}
\lim_{N\rightarrow\infty} \mathcal{R}_t^N(\al)=
\frac{\Omega_t(\al)}{(\Omega_0(\al))^{1-t}(\Omega_1(\al)^t},
\end{equation}
uniformly for  $\alpha \in [0, 1]$.
\end{corollary}

The next Corollary follows immediately from Theorem \ref{main1}
and Lemma \ref{LBCQ}.

\begin{corollary} \label{cora}

There exist $C>0$ such that for all $\alpha \in [0,1]\cap
\frac{1}{N}\mathbf{Z}$

\begin{equation}\label{polybound1}
Q_t^N(\alpha)\geq C.
\end{equation}

\end{corollary}

We divide the proof of Theorem \ref{main1} into an analysis of
norming constants in an interior region of $[0, 1]$ and in a
boundary region.


\subsection{Interior estimates}

We begin by studying $Q_t^N(\alpha)$ where $\alpha$ lies in the (left) `interior interval' $\alpha \in
[\frac{1}{N^{3/4}}, \frac{2}{3}]$. It is then possible to obtain
joint $(N, \alpha)$ asymptotics by a complex stationary phase
method. The discussion is essentially the same for the right
interior interval $[\frac{1}{3}, 1 - \frac{1}{N^{3/4}}]$ and is
omitted.

\begin{proposition}\label{asy}

 Let $\alpha \in [\frac{1}{N^{3/4}}, \frac{2}{3}]$ and $M=N\alpha$. Then there exist uniformly bounded functions $A_{t, k}(\alpha)$
 on the interior region, such that
\begin{equation}
Q_t^N(\alpha) \sim \frac{2 \pi ^{\frac{3}{2}}\alpha}{(\frac{1}{(1-\alpha)}+ \alpha
f''_t(\alpha))^{\frac{1}{2}}(M)^{\frac{1}{2}}}\sum_{k=0}^{\infty}
\frac{A_{t,k}(\alpha)}{M^k}=
\frac{2 \pi ^{\frac{3}{2}}}{(G_t''( \alpha))^{\frac{1}{2}}(N)^{\frac{1}{2}}}\sum_{k=0}^{\infty}
\frac{A_{t,k}(\alpha)}{M^k},
\end{equation}
in the sense that for any $R\in \mathbf{Z}^+$  there exists
$C_R>0$ such that
\begin{equation}
\left|Q_t^N(\alpha) -\frac{2 \pi ^{\frac{3}{2}}\alpha}{(\frac{1}{(1-\alpha)}+ \alpha
f''_t(\alpha))^{\frac{1}{2}}(M)^{\frac{1}{2}}}\sum_{k=0}^{R}
\frac{A_{t,k}(\alpha)}{M^k}\right|\leq
\frac{C_R \alpha }{(\frac{1}{(1-\alpha)}+ \alpha
f''_t(\alpha))^{\frac{1}{2}}(M)^{\frac{1}{2}}} M^{-(R+1)}.
\end{equation}
In particular, $A_{t,0}=1$.

\end{proposition}

\begin{corollary} \label{inasy} Let $\alpha \in [\frac{1}{N^{3/4}}, \frac{2}{3}]$ and $M=N\alpha$. There is a complete asymptotic expansion for large $M$

\begin{equation}
q_t^N(\alpha)\sim \frac{1}{(M
(1+\alpha(1-\alpha)f''_t(\alpha)))^{\frac{1}{2}}}
\sum_{k=0}^{\infty} \frac{B_{t,k}(\alpha)}{M^k}=
\left(\frac{G_{\PP}''(\alpha)}{G''_t(\alpha)}\right)^{\frac{1}{2}}
\sum_{k=0}^{\infty} \frac{B_{t,k}(\alpha)}{M^k}
\end{equation}
in the sense that for any $R\in
\mathbf{Z}^+$ there exists $C_R>0$ such that
\begin{equation}
\left|q_t^N(\alpha) -\frac{1}{(M (1+\alpha(1-\alpha)
f''_t(\alpha)))^{\frac{1}{2}}}\sum_{k=0}^{R}
\frac{B_{t,k}(\alpha)}{M^k}\right|\leq C_R M^{-(R+1)}.
\end{equation}
In particular, $B_{t,0}=1$ and there exists $C>0$ such that

\begin{equation}
0< \frac{1}{C} \leq q_t^N(\alpha)\leq C.
\end{equation}

\end{corollary}

The proof of  Proposition \ref{asy}
 proceeds by a sequence of Lemmas. The first concerns the phase
 $F_{t, \alpha}$ (\ref{FTALPHA}).

\begin{lemma}\label{bound2}
$\alpha$ is the only critical point of $F_{t,\alpha}(x)$ and we have
\begin{equation}
F''_{t,\alpha}(\alpha)=G_t''(\alpha)>0, ~~~~
(x-\alpha)F'_{t,\alpha} (x)\geq0.
\end{equation}

\end{lemma}

\begin{proof} Differentiating (\ref{PHASEDEF}) shows that  $F'_{t,\alpha}(x)=(x-\alpha)
G''_t(x)$. The second derivative is readily obtained and it is
positive by Lemma \ref{bound1}.

\qed

\end{proof}

Now we make a substitution of variables. Let
$y=\frac{x-\alpha}{\alpha}$, $M=N\alpha$. We then have
\begin{eqnarray} \label{QTALPHA} 
Q_t^N(\alpha) &=& 2\pi \alpha \int_{-1}^{\frac{1}{\alpha}-1}e^{-M
\FF_{t,\alpha}(y)}dy
\end{eqnarray}
with new phase function

\begin{equation}
\begin{array}{lll} \FF_{t, \alpha}(y) & = &
\frac{1}{\alpha}F_{t,\alpha}(\alpha(1+y)), ~~ \FF_{P, \alpha}(y) =
\frac{1}{\alpha}F_{P,\alpha}(\alpha(1+y)) \\ && \\
\FF_{t, \alpha}(y)&=&  -\left(\log
(1+y)+\frac{1-\alpha}{\alpha}\log\frac{1-\alpha-\alpha
y}{1-\alpha}+\alpha y^2 f_{t,\alpha}(\alpha(1+y))\right).
\end{array} \end{equation}

\begin{lemma}\label{newphase} The phase has the following properties:
\begin{enumerate}

\item  $\FF_{t,\alpha}(y)$ is strictly decreasing on $(-1, 0)$ and
strictly increasing on $(0, \frac{1}{\alpha}-1)$ with a unique
critical (minimum) point  at  $y=0$ with $\FF_{t,\alpha}(0) = 0$.

\item If $y_0 > 0$, then   $\inf_{y \geq y_0} \FF'_{t,\alpha}(y)
\geq C(y_0) > 0$ where $C(y_0)$ is independent of $\alpha, t$.

\item If $y_0 < 0$, then $\inf_{y \in [-1, y_0]}
|\FF'_{t,\alpha}(y)| \geq C(y_0) > 0$ where $C(y_0)$ is
independent of $\alpha, t$.

\item The Hessian of $\FF_{t, \alpha}$ of $y = 0$ is
non-degenerate and
$$
\FF''_{t,\alpha}(0) =\alpha G_t''(\alpha) = \frac{1}{1-\alpha} +
\alpha f_t''(\alpha) >0. $$

\item  $\FF_{t,\alpha}(y)$ and all of
 its derivatives are uniformly bounded for $\alpha\in [0,
\frac{2}{3}]$ and for $y$ in   any compact set of $(-1,
\frac{1}{\alpha}-1)$.

\end{enumerate}

\end{lemma}

\begin{proof} Comparing with (\ref{FTALPHA}) and Lemma \ref{bound2} shows that
\begin{eqnarray*}\label{trythislabel}
\FF_{t,\alpha}'(y)&=& \frac{1}{\alpha}\frac{d F_{t,
\alpha}(x)}{dx}\frac{dx}{dy} =F_{t, \alpha} '(x) \\ %
 &&  =
(x-\alpha) G_t''(x) = \alpha y G_t''(\alpha(1+y))
=\frac{x-\alpha}{x}(xG''_{t}(x))=\frac{y}{1+y}
(xG''_t(x)) \\ %
 & = & -\frac{1}{1+y}
+\frac{1-\alpha}{1-\alpha-\alpha y}+2\alpha y
f_{t, \alpha}(\alpha(1+y))+\alpha^2 y^2 f'_{t,\alpha}(\alpha(1+y))\\
\FF_{t,\alpha}''(y)&=& \alpha F''_{t, \alpha} (x)=\alpha G_t''(x) + \alpha (x-\alpha)G_t'''(x) =\alpha G_t''(\alpha(1+y)) + \alpha^2 y G_t'''(\alpha(1+y))\\
&=&\frac{1}{(1+y)^2}
+\frac{\alpha(1-\alpha)}{(1-\alpha-\alpha y)^2}\\
&&+2 \alpha f_{t,\alpha}(\alpha(1+y))+ (2 \alpha y+2\alpha^2 y)
f'_{t,\alpha}(\alpha(1+y))+\alpha^2y^2 f_{t,\alpha}(\alpha(1+y)).
\end{eqnarray*}

By Lemma \ref{bound1}, $xG''_t(x)$ has a uniform positive lower
bound, hence by the formula $\FF_{t,\alpha}'(y) = \frac{y}{1+y}
(xG''_t(x))$,  $\FF_{t,\alpha}'(y)=0$ if and only if $y=0$. Also
$\FF_{t,\alpha}'(y)<0$ on $(-1,0)$ and $\FF_{t,\alpha}'(y)>0$ on
$( 0 , \frac{1}{\alpha}-1 )$. The same formula implies (2)-(3)
since the factor $|\frac{y}{1+y}|$ then has a uniform lower bound.

Again by Lemma \ref{bound1},  $G_t''(\alpha)$ has poles at $0$ and
$1$, hence $\alpha G_t''(\alpha)$ is uniformly bounded below from
$0$ for $\alpha\in [0, \frac{2}{3}]$. In particular, at the
critical point, we have (cf. Lemma \ref{bound2}),
\begin{eqnarray*}
\FF''_{t,\alpha}(0) = \alpha F_{t, \alpha}''(\alpha)  =
\frac{1}{1-\alpha} + \alpha f_t''(\alpha)=\alpha G_t''(\alpha)>0.
\end{eqnarray*}

\qed
\end{proof}

\begin{lemma}\label{local1}

There exist $\delta$ and $C>0$ such that
\begin{equation}
\left|1-\frac{2\pi \alpha\int_{-\frac{1}{2}}^{1}e^{-M
\FF_{t,\alpha}(y)}dy}{Q_t^N(\alpha)} \right| \leq  \frac{ C e^{-\delta
M}}{M}.
\end{equation}

\end{lemma}

\begin{proof}

By Lemma \ref{newphase} (2),  there exists   $\Lambda>0$
independent of $(t, \alpha)$ such that
\begin{equation}\label{F'}\FF_{t,\alpha}(y)  \geq \FF_{t,\alpha}(1) + \frac{\Lambda}{2} (y - 1),\;\;\mbox{for}\;  y \geq 1. \end{equation}
 Using also
 that $\FF_{t,\alpha}$ increases on $(0, \frac{1}{2})$, we have
\begin{eqnarray*}\nonumber
\int_1^{\frac{1}{\alpha}-1}e^{-M\FF_{t,\alpha}(y)}dy &\leq&
\int_{1}^{\frac{1}{\alpha}-1}
e^{-\frac{\Lambda}{2} M (y-1)-M\FF_{t,\alpha}(1)}dy\\
&\leq&\frac{2e^{-M\FF_{t, \alpha}(1)}}{\Lambda M}\\
&\leq& \frac{4e^{-M(\FF_{t, \alpha}(1)-\FF_{t,
\alpha}(\frac{1}{2}))}}{\Lambda
M}\int_{0}^{\frac{1}{2}}e^{-M\FF_{t,\alpha}(y)}dy\\
&\leq& \frac{Ce^{-\delta M}}{2M\alpha} Q_t^N(\alpha),\;\;
\mbox{where}\;\; \delta: = 2 \inf_{y \in [\frac{1}{2}, 1]}
\FF_{t, \alpha}'(y).
\end{eqnarray*}
In the last line we again used Lemma \ref{newphase} (2).

By the same argument, there exists $\delta > 0$ (independent of
$(t, \alpha)$ so that
$$\int_{-1}^{-\frac{1}{2}}e^{-M\FF_{t,\alpha}(y)}dy
\leq \frac{Ce^{-\delta M}}{2M\alpha} Q_t^N(\alpha).$$ Indeed, by
Lemma \ref{bound1} (3), $\FF_{t,\alpha}$ is decreasing on $(-1,
0)$ and there exists $ - \Lambda < 0$ independent of $(t, \alpha)$
so that
\begin{equation}\label{F'2} \FF_{t,\alpha}(y)  \geq \FF_{t,\alpha}(-\frac{1}{2}) - \frac{\Lambda}{2} (y + \frac{1}{2}). \end{equation}
As above,
\begin{eqnarray*}\nonumber
\int_{-1}^{-\frac{1}{2}} e^{-M\FF_{t,\alpha}(y)}dy
&\leq&\frac{2e^{-M\FF_{t, \alpha}(-\frac{1}{2})}}{\Lambda M}\\
&\leq& \frac{8 e^{-M(\FF_{t, \alpha}(-\frac{1}{2})-\FF_{t,
\alpha}(-\frac{1}{4}))}}{\Lambda
M}\int_{-\frac{1}{2}}^{-\frac{1}{4}}e^{-M\FF_{t,\alpha}(y)}dy\\
&\leq& \frac{Ce^{-\delta M}}{2M\alpha} Q_t^N(\alpha),\;\;
\mbox{where}\;\; \delta: = 2 \inf_{y \in [-\frac{1}{2}, -
\frac{1}{4}]} |\FF_{t, \alpha}'(y)|.
\end{eqnarray*}

The lemma is proved by combining the above inequalities and

$$\left|1-\frac{2\pi \alpha\int_{-\frac{1}{2}}^{1}e^{-M
\FF_{t,\alpha}(y)}dy}{Q_t^N(\alpha)} \right| =  \left|
 \frac{ 2\pi \alpha\int_{-1}^{-\frac{1}{2}}e^{-M
\FF_{t,\alpha}(y)}dy + 2\pi \alpha\int_{1}^{\frac{1}{\alpha} - 1}e^{-M
\FF_{t,\alpha}(y)}dy}{Q_t^N(\alpha)}
 \right|   \leq  \frac{ C e^{-\delta
M}}{M}.$$

\qed
\end{proof}

\noindent{\bf Proof of Proposition \ref{asy}}\\

\noindent

We introduce a smooth cut-off function $\eta$ such that $\eta=1$
on $[-\frac{1}{2} + \epsilon, 1 - \epsilon]$ for some fixed
$\epsilon > 0$ (independent of $(t, N, \alpha))$  and with
$\eta=0$ outside $(-\frac{1}{2}, 1)$ and write
$$Q_t^N(\alpha) = I_t^N(\alpha) + II_t^N(\alpha), \;\; \mbox{with}\;\; I_t^N(\alpha) = 2\pi \alpha  \int_{-\frac{1}{2}}^{1} e^{-M \FF_{t, \alpha}(y)}
\eta(y)  dy. $$

By Lemma \ref{local1}, $II_t^N \leq  \frac{Ce^{-\delta
M}}{2M\alpha} Q_t^N(\alpha), $  hence $Q_t^N(\alpha) \left(1 +
O(\frac{e^{-\delta M}}{2M\alpha})\right) = I_t^N(\alpha), $ and
therefore \begin{equation} \label{QI}  Q_t^N(\alpha)
  =I_t^N(\alpha) \left(1 + O(\frac{e^{-\delta M}}{2M\alpha})\right).
\end{equation}

We now evaluate  $I_t^N(\alpha)$ asymptotically with respect to
the parameter $M$ by  the method of complex stationary phase with
non-degenerate complex phase functions (Theorem 7.7.5 in
\cite{H}), with $\alpha, t$ as parameters. Using the evaluation of
the Hessian in  Lemma \ref{newphase},  we obtain
\begin{equation}
\left|I_t^N(\alpha) -\frac{2 \pi ^{\frac{3}{2}}\alpha}{(\frac{1}{(1-\alpha)}+ \alpha
f''_t(\alpha))^{\frac{1}{2}}(M)^{\frac{1}{2}}}\sum_{k=0}^{R}
\frac{A_{t,k}(\alpha)}{M^k}\right|\leq \frac{C_R \alpha
}{(\frac{1}{(1-\alpha)}+ \alpha
f''_t(\alpha))^{\frac{1}{2}}(M)^{\frac{1}{2}}} M^{-(R+1)},
\end{equation}
where the $A_{t, k}(\alpha)$  are obtained by applying powers of
the inverse Hessian operator
$$\frac{1}{\left(\frac{1}{(1-\alpha)}+
\alpha f''_t(\alpha) \right)} \frac{d^2}{dy^2} $$ to $e^{-M R_3(y;
t, \alpha)}$ and $R_3$ is the cubic remainder in the Taylor
expansion of $\FF_{t, \alpha}(y)$ at $y = 0$. The inverse Hessian
is uniformly bounded above in the interior region, so $R_3$ is uniformly bounded with
uniformly bounded derivatives when $\alpha \in [0, \frac{2}{3}]$
and $y \in [-\frac{1}{2} + \epsilon, 1 - \epsilon]$ (cf. Lemma
\ref{newphase} (5). Therefore the stationary phase
coefficients and remainder are uniformly bounded in the interior
region.

The Proposition follows by combining the complex stationary phase
asymptotics with (\ref{QI}).

\qed

\subsection{Boundary estimates}

We now give joint asymptotic estimates of norming constants in the
boundary zone where $0< \alpha \leq \frac{1}{N^{3/4}}$. The
exclusion of $\alpha = 0$ is not important since the norming
constants also equal $1$ there. The main result of this section
is:

\begin{proposition}\label{bdasy} Assume $0<\alpha \leq \frac{1}{N^{3/4}}$. Then we have
\begin{equation}
q_t^N(\alpha)=1+O(N^{-\frac{1}{3}}).
\end{equation}
\end{proposition}

The proof of Proposition \ref{bdasy} consists of a number of
Lemmas.

First we will localize the integral $Q_t^N$ and $Q_P^N$.
\begin{lemma} \label{bdylem} Suppose $0< \alpha \leq \frac{1}{N^{3/4}}$. Then there
exists constants $\delta, C>0$ such that

\begin{eqnarray}
\left\{\begin{array}{lll}
\frac{\alpha\int_{\frac{1}{\alpha N^{2/3}}}^{\frac{1}{\alpha}-1}
e^{-\alpha N \FF_{t,\alpha}(y)}dy}{Q_t^N(\alpha)} &\leq& C
e^{-\delta N^{1/3}},\\
&&\\
\frac{\alpha\int_{\frac{1}{\alpha N^{2/3}}}^{\frac{1}{\alpha}-1}
e^{-\alpha N \FF_{P,\alpha}(y)}dy}{Q_P^N(\alpha)} &\leq& C
e^{-\delta N^{1/3}}.
\end{array} \right.
\end{eqnarray}
\end{lemma}

\begin{proof} We now localize   the integral (\ref{QTALPHA}) into
the subinterval $[-1, \frac{1}{\alpha N^{2/3}}]$ by showing that
the integral over $[\frac{1}{\alpha N^{2/3}}, \frac{1}{\alpha} -
1]$ is relatively negligeable. In the boundary region,
$\frac{1}{\alpha N^{2/3}} \geq N^{1/12}$.

As in  Lemma \ref{local1}, there exists a uniform positive
constant $\Lambda > 0$  such that $\FF_{t,\alpha}(y) \geq
\FF_{t,\alpha}(\frac{1}{\alpha N^{2/3}}) + \Lambda (y -
\frac{1}{\alpha N^{2/3}})$ on  $[\frac{1}{\alpha N^{2/3}},
\frac{1}{\alpha} - 1]$ and we obtain
\begin{eqnarray*}
\int_{\frac{1}{\alpha N^{2/3}}}^{\frac{1}{\alpha}-1} e^{-\alpha N
\FF_{t,\alpha}(y)}dy &\leq&\frac{C}{\alpha N}e^{-\alpha N
\FF_{t,\alpha}(\frac{1}{\alpha N^{2/3}})}\\
&\leq&\frac{C}{\alpha N}e^{-\epsilon N^{1/3} },
\end{eqnarray*}
using the fact (cf. Lemma \ref{newphase}) that $\FF_{t,\alpha}(y)
\geq
 \epsilon (y-1)  + \FF_{t,\alpha}(1)$ as $y \to \infty$; here, $\epsilon > 0$  is a positive and uniform constant.

To prove that the integral over $[\frac{1}{\alpha N^{2/3}},
\frac{1}{\alpha} - 1]$ is relatively negligeable, we   give a
lower bound for the integral on the whole interval $[-1,
\frac{1}{\alpha} - 1]$. In fact a very crude lower bound suffices,
so we choose a convenient subinterval. By Lemma \ref{newphase},
$\FF_{t,\alpha}'(y)$ is uniformly bounded on any compact subset of
$[0,\frac{1}{\alpha}-1]$ when $\alpha$ is the boundary region.
Using that $\FF_{t, \alpha}(y) \leq \FF_{t, \alpha}(1)  + C_0 (y-1)$ on $[1,2]$ for some $C_0>0$, 
we have
\begin{eqnarray*}
\int_{-1}^{\frac{1}{\alpha}-1} e^{-\alpha N \FF_{t,\alpha}(y)}dy
&\geq & \int_{1}^{2} e^{-\alpha N \FF_{t,\alpha}(y)}dy \\
&\geq& e^{-\alpha N
\FF_{t,\alpha}(1)}\int_{1}^{2}e^{-C_0\alpha N (y-1)}dy\\
&\geq&\frac{C}{\alpha N}e^{-\alpha N \FF_{t,\alpha}(1)-C_0\alpha
N}\\
&\geq&\frac{C}{\alpha N}e^{-C_0 N^{1/4}}.
\end{eqnarray*}

\noindent Therefore
\begin{eqnarray*}
\frac{\alpha\int_{\frac{1}{\alpha N^{2/3}}}^{\frac{1}{\alpha}-1}
e^{-\alpha N \FF_{t,\alpha}(y)}dy}{Q_t^N(\alpha)}\leq C
e^{-\epsilon N^{1/3}+C_0 N^{1/4}}\leq C e^{- \delta N^{1/3}}
\end{eqnarray*}
for some $\delta>0$.

\qed
\end{proof}

\subsection{Proof of Proposition \ref{bdasy} }

\begin{proof} By definition and the previous lemma, we have

\begin{eqnarray*}
&&q_t^N(\alpha)=\frac{Q_t^N(\alpha)}{Q_P^N(\alpha)}\\
&=&\frac{\int^{\frac{1}{\alpha N^{2/3}}}_{-1} e^{-\alpha N
\FF_{t,\alpha}(y)}dy}{\int^{\frac{1}{\alpha N^{2/3}}}_{-1}
e^{-\alpha N \FF_{P,\alpha}(y)}dy}
(1+O(e^{N^{-\frac{1}{3}}}))\\
&=&\frac{\int^{\frac{1}{\alpha N^{2/3}}}_{-1} e^{-\alpha N
\FF_{P,\alpha}(y)+\alpha^2N y^2
f_{t,\alpha}(\alpha(1+y))}dy}{\int^{\frac{1}{\alpha N^{2/3}}}_{-1}
e^{-\alpha N \FF_{P,\alpha}(y)}dy}
(1+O(e^{- N^{-\frac{1}{3}}}))\\
&=&(1+O(N^{-\frac{1}{3}})),
\end{eqnarray*}
where in the last line, we Taylor expand the exponential
$e^{-\alpha^2N y^2 f_{t,\alpha}(\alpha(1+y))}$ and recall that in
the boundary layer, $\alpha^2N y^2  f_{t,\alpha}(\alpha(1+y)) =
O(N^{-\frac{1}{3}})$.

This completes the proof of Proposition
\ref{bdasy}.

\end{proof}

When $\alpha=0$,

$$ q_t^N(0)= \frac{ \int_{0}^{1} e^{N\log(1-x)-Nx^2 f_{t,\alpha}(x)}dx}
{\int_{0}^{1} e^{N\log(1-x)}dx}. $$
Notice that the phase is strictly decreasing on $[0,1]$, one can apply the similar argument as before with the substitution $y=N^{2/3}x$. We leave it as an excercise to show
$q_t^N(0)\sim 1$ as $N\to \infty$.

\qed

\subsection{\bf Completion of proof of Theorem \ref{main1}}

 It is easy to see in the
boundary layer $0< \alpha \leq \frac{1}{N^{3/4}}$,
$$\frac{\det\nabla^2G_t(\alpha)}
{\det\nabla^2 G_P(\alpha)}=
\frac{\frac{1}{\alpha(1-\alpha)}+f''_t(\alpha)}{\frac{1}{\alpha(1-\alpha)}}
=1+O(\alpha).$$
Therefore $\frac{Q_t^N(\alpha)}{Q_P^N(\alpha)}$ continuously extends to $P$.

 Consider the interior of $P$ with $\frac{1}{N^{3/4}}\leq
\alpha\leq \frac{2}{3}$. By  Corollary \ref{inasy}, we have

$$q_t^N(\alpha)=\frac{Q_t^N(\alpha)}{Q_P^N(\alpha)}=\left(\frac{1}{1+f''_t(\alpha)}\right)^
{\frac{1}{2}}+O(\frac{1}{\alpha N})=\frac{\det\nabla^2G_t(\alpha)}
{\det\nabla^2 G_P(\alpha)}+O(\frac{1}{\alpha
N})=\frac{\det\nabla^2G_t(\alpha)} {\det\nabla^2
G_P(\alpha)}+O(\frac{1}{N^{1/4}}).$$

 Consider the boundary layer of $P$ with
$0\leq \alpha\leq\frac{1}{N^{3/4}}$. By Lemma \ref{bdasy}, we have
$$q_t^N(\alpha)=\frac{Q_t^N(\alpha)}{Q_P^N(\alpha)}=1+O(\frac{1}{N^{1/3}})=\frac{\det\nabla^2G_t(\alpha)} {\det\nabla^2
G_P(\alpha)}+O(\frac{1}{N^{1/3}}).$$

\noindent  Therefore for any $\alpha \in P$, we have
\begin{equation}
q_t^N(\alpha)=\frac{\det\nabla^2G_t(\alpha)} {\det\nabla^2
G_P(\alpha)}+O(\frac{1}{N^{1/4}}).
\end{equation}
This proves Theorem \ref{main1}.

\qed


\section{Proof of the main theorem \ref{main}}

\subsection{The $C^0$-convergence}

\begin{proposition}\label{C0}  There exists a constant $C>0$ such that for
any $z\in \CPP^1$

\begin{equation}
\left|\frac{1}{N}\log \CE_N(t, z)\right|\leq \frac{C\log N}{N}.
\end{equation}

\end{proposition}

\begin{proof}

There exists a constant $C>0$ such that $\frac{1}{C}\leq
\CR_t^N(\alpha) \leq C$ for any $\alpha\in P$. Then
$$ \frac{1}{C}\leq\frac{\CE_N(t, z)}{\Pi_N(t, z)}\leq C.$$
The proposition is proved by applying the Tian-Yau-Zelditch
expansion \cite{Z}, which asserts that there exists a $C^{\infty}$
asymptotic expansion,
\begin{equation} \label{TYZ} \Pi_N(t, z)= N \; \left(1 +
\frac{a_1(t, z)}{N} + \cdots \right),
\end{equation} which may be differentiated any number of times. It obviously implies that
$$ \left| \frac{1}{N} \log \CE_N(t, z) - \frac{1}{N} \log \Pi_N(t,z) \right| \leq \frac{ C \log N}{N}. $$

\qed
\end{proof}

We note that the convergence rate is best possible. The $C^0$
convergence with this rate was proved for general \kahler
manifolds in \cite{B}; a simple proof along the above lines  for
general toric \kahler manifolds will appear in \cite{SoZ}.

\subsection{\bf A localization lemma}

To obtain $C^2$ convergence, we have to estimate weighted sums of
$\CP_{\alpha}^N(t, z)$ for $\alpha \in P\cap \frac{1}{N}
\mathbf{Z}$. The following localization lemma enables to replace
the sum of $\CP_{\alpha}^N(t, z)$ by its partial sum for $\alpha$
in  small neighborhoods of $\pi_t( \rho )$, where $\rho=\log
|z|^2$. Fix $ z = e^{\rho/2+i\theta }\in X$ with $x=\pi_t (\rho)$.

\begin{lemma}\label{loc}

 For any $\delta>0$, there exist  $0<\delta_1<\delta$, $0<\delta_2<\delta$, $\epsilon>0$ and $C>0$ such that for any $\alpha$ and $\beta \in [0,1] \cap \frac{1}{N} \mathbf{Z}$ with $ |\alpha-x| < \delta_1$ and $|\beta-x|\geq 2 \delta_2$, we have

\begin{equation}
\frac{\PP_{\alpha}^N (t,z)}{\PP_{\beta}^N(t,z)}\leq C e^{-\epsilon N}.
\end{equation}

\end{lemma}

\begin{proof} First let's assume $x\in (0,1)$.

\begin{eqnarray*}
\frac{ \PP_{\beta}^N(t,z)}{ \PP_{\alpha}^N (t, z) } &=& \frac{ e^{-N \left( (x-\beta) G_t'(x)-(G_t(x)-G_t(\beta))\right)} } { e^{-N \left( (x-\alpha) G_t'(x)-(G_t(x)-G_t(\alpha))\right) }}
\frac{ Q_t^N(\alpha)} {Q_t^N(\beta)}\\
&=& e^{ -N  \left(  G_t(\beta)-G_t(\alpha)- G'_t(x)(\beta -  \alpha) \right) }\frac{ Q_t^N(\alpha)} {Q_t^N(\beta)}\\
&=& e^{ -N \left( G_t(\beta) - G_t (\alpha) - G'_t(\alpha) (\beta-\alpha) \right) + N(G'_t(x)-G'_t(\alpha))(\beta-\alpha)} \frac{ Q_t^N(\alpha)} {Q_t^N(\beta)}\\
&=&e^{ -N G''_t(\gamma) (\beta-\alpha)^2  + N(G'_t(x)-G'_t(\alpha))(\beta-\alpha)} \frac{ Q_t^N(\alpha)} {Q_t^N(\beta)}
\end{eqnarray*}
for some $\gamma$ between $\alpha$ and $\beta$.

Notice that $G''_t$ is uniformly bounded below from $0$ and $G'_t$ is equicontinuous on [0,1]. Therefore we can choose $\delta_1 << \delta_2$ so that there exits $\epsilon>0$ and $$ e^{ -N G''_t(\gamma) (\beta-\alpha)^2  + N(G'_t(x)-G'_t(\alpha))(\beta-\alpha)} \leq e^{-2\epsilon N}.$$

Also $$  \frac{\PP_{\alpha}^N (t,z)}{\PP_{\beta}^N(t,z)}\leq C e^{-\epsilon N}. $$

When $x=0$ or $1$,  the same estimate can be proved by similar argument as above making use of the monotonicity of $G_t$.

\qed
\end{proof}


\bigskip

\subsection{The $C^{2}$-convergence}

We now prove the main results giving bounds on two space-time
derivatives of $\varphi_N(t, z)$. The main ingredients are the
bounds of $\CR_t^N$ in Lemma \ref{maincor1} and a comparison to
derivatives of the \szego kernel $\Pi_N(t, z)$ for the metric
$\phi_t$ deriving from Lemma \ref{RCAL}. By (\ref{TYZ}) it is
straightforward to determine the derivatives of $\Pi_N(t,z)$.

The following lemma is the consequence of the family version of the Tian-Yau-Zelditch expansion.

\begin{lemma} We have the following uniform convergence in the $C^{\infty}$ topology on $[0,1]\times \CPP^1$
\begin{equation}
\lim_{N\rightarrow \infty} \frac{1}{N} \log \Pi_N(t, z)=0.
\end{equation}

\end{lemma}

\begin{corollary}  All derivatives of $   \frac{1}{N} \log \Pi_N (t, z) + u_t(z)  $, of order great than zero,  are uniformly bounded on $X$.

\end{corollary}

\begin{proof}

Although $ u_t(z)  $ is not a well-defined function on $X$, $e^{-u_t(z)}$ extends to a hermitian metric on the line bundle so that, by applying global vector fields, any derivatives of $u_t(z)$ are well defined functions on $X$ and are uniformly bounded.

\qed
\end{proof}

\begin{proposition} \label{C2BOUNDS}

\begin{equation}
\lim_{N\rightarrow\infty} \left|\left|\frac{1}{N}\log \CE_N(t,
z)\right|\right|_{C^{2} ([0,1]\times X)}=0.
\end{equation}

\end{proposition}

\begin{proof}   Fix $z\in \CPP^1$, and put $x=\pi_t(z)$. To prove
the $C^2$ convergence of $\frac{1}{N} \log \ecal_N(t,z)$ it
suffices by (\ref{TYZ}) to prove $C^2$ convergence for
$\frac{1}{N} ( \log \ecal_N(t,z)- \log \Pi_N(t, z))$.  We use
Lemma \ref{RCAL} to simplify the formula for $\CE_N(t,z)$.

\bigskip

\noindent {\bf Second order convergence in pure space derivatives}

\bigskip

We first consider pure space derivatives.  By $\sum_{\alpha}$ and $\sum_{\alpha, \beta}$, we mean $\sum_{\alpha \in P\cap \frac{1}{N}\mathbf{Z}}$ and $\sum_{\alpha, \beta \in P\cap \frac{1}{N}\mathbf{Z}}$. If $x$ is in the
`interior region' of $[0,1]$, we may use the coordinates $z =
e^{\rho/2 + i \theta}$, and
\begin{eqnarray*}
&&\frac{1}{N}\left|\frac{\partial^2}{\partial \rho^2} \log
\ecal_N(t,z)-\frac{\partial^2}{\partial \rho^2} \log \Pi_N(t, z)\right|\\
&=&N\left|\frac{\sum_{\alpha,\beta} (\alpha-\beta)^2
\CR_t^N(\alpha)\CR_t^N(\beta) \CP^N_{\alpha}(t,z)\CP^N_{\beta}(t,z)}{2\left(\sum_{\alpha} \CR_t^N(\alpha)\CP^N_{\alpha}(t,z)\right)^2} 
-\frac{\sum_{\alpha,\beta} (\alpha-\beta)^2
 \CP^N_{\alpha}(t,z)\CP^N_{\beta}(t,z)}{2\left(\sum_{\alpha} \CP^N_{\alpha}(t,z))\right)^2} \right|\\
&=&N\left(\frac{\sum_{\alpha,\beta} (\alpha-\beta)^2
\CP^N_{\alpha}(t,z)\CP^N_{\beta}(t,z)}{2\left(\sum_{\alpha} \CP^N_{\alpha}(t,z)\right)^2}\right) \\
&&\left|\frac{\sum_{\alpha,\beta} (\alpha-\beta)^2
\CR_t^N(\alpha)\CR_t^N(\beta) \CP^N_{\alpha}(t,z)\CP^N_{\beta}(t,z)}{\left(\sum_{\alpha} \CR_t^N(\alpha)\CP^N_{\alpha}(t,z)\right)^2}  \frac{\left(\sum_{\alpha} \CP^N_{\alpha}(t,z)\right)^2}{\sum_{\alpha,\beta} (\alpha-\beta)^2
\CP^N_{\alpha}(t,z)\CP^N_{\beta}(t,z)} -1\right| \\
&\leq& \left|\frac{1}{N} \frac{\partial^2}{\partial \rho^2} \log
\Pi_N(t,z) + \frac{\partial^2}{\partial \rho^2} u_t(z)\right| \\
&&\left|\frac{\sum_{\alpha,\beta} (\alpha-\beta)^2
\CR_t^N(\alpha)\CR_t^N(\beta)\CP^N_{\alpha}(t,z)\CP^N_{\beta}(t,z)}{\left(\sum_{\alpha} \CR_t^N(\alpha)\CP^N_{\alpha}(t,z)\right)^2}  \frac{\left(\sum_{\alpha} \CP^N_{\alpha}(t,z)\right)^2}{\sum_{\alpha,\beta} (\alpha-\beta)^2
\CP^N_{\alpha}(t,z)\CP^N_{\beta}(t,z)} -1\right| \\
&\leq & C\left|\frac{\sum_{\alpha,\beta}
(\alpha-\beta)^2 \CR_t^N(\alpha)\CR_t^N(\beta) \CP^N_{\alpha}(t,z)\CP^N_{\beta}(t,z)}{\sum_{\alpha,\beta} (\alpha-\beta)^2
\CP^N_{\alpha}(t,z)\CP^N_{\beta}(t,z)}
\frac{\left(\sum_{\alpha} \CP_t^N(\alpha,
z)\right)^2}{\left(\sum_{\alpha}
\CR_t^N(\alpha)\CP^N_{\alpha}(t,z)\right)^2} -1\right| \\
&\leq & C\left|\frac{\sum_{\alpha,\beta \in B_x(\delta)}
(\alpha-\beta)^2 \CR_t^N(\alpha)\CR_t^N(\beta) \CP^N_{\alpha}(t,z)\CP^N_{\beta}(t,z)}{\sum_{\alpha,\beta \in B_x(\delta)} (\alpha-\beta)^2
\CP^N_{\alpha}(t,z)\CP^N_{\beta}(t,z)}
\frac{\left(\sum_{\alpha \in B_x(\delta)} \CP_t^N(\alpha,
z)\right)^2}{\left(\sum_{\alpha \in B_x(\delta)}
\CR_t^N(\alpha)\CP^N_{\alpha}(t,z)\right)^2} -1 + O\left( e^{-\epsilon N} \right)\right| \\
&\leq& C \left( \sup_{\alpha\in B_x(\delta)} \CR_t^N(\alpha) -\inf_{\alpha\in B_x(\delta)} \CR_t^N(\alpha) \right) +O\left( e^{-\epsilon N} \right)
\end{eqnarray*}
for some fixed $\epsilon>0$,  where $B_x(\delta)=\{ \alpha\in [0,1]\cap \frac{1}{N}\mathbf{Z} ~|~|\alpha- x|<\delta \}$.

$R_t^{\infty}(\alpha) = \lim_{N\rightarrow \infty} R^N_t(\alpha)$ is continuous on $[0,1]$ and the convergence is uniform.
 Therefore for any $\epsilon' >0$, there exists $\delta>0$ and sufficiently large $N'$ such that for all $N\geq N'$
$$  \sup_{\alpha\in B_x(\delta)} \CR_t^N(\alpha) - \inf_{\alpha\in B_x(\delta)} \CR_t^N(\alpha) \leq \epsilon'.$$

In other words, for any $\epsilon' >0$, there exists a sufficiently large $N'$ such that for all $N\geq N'$,
$$\frac{1}{N}\left|\frac{\partial^2}{\partial \rho^2} \log
\ecal_N(t,z)-\frac{\partial^2}{\partial \rho^2} \log \Pi_N(t, z)\right|< \epsilon'.$$

If $x$ is close to the boundary of $[0,1]$, without loss of generality  we fix a holomorphic coordinate system $\{z\}$ for $\mathbf{CP}^1$ near the north
pole $z_N$ such that $z=0$ at $z_N$. Let $r=|z|$. Then

\begin{eqnarray*}
&&\frac{1}{N}\left|\frac{\partial^2}{\partial r^2} \log
\ecal_N(t,z)-\frac{\partial^2}{\partial r^2} \log \Pi_N(t,z)\right|\\
&=&N\left|\frac{\sum_{\alpha,\beta>0} (\alpha-\beta)^2
\CR_t^N(\alpha)\CR_t^N(\beta) \CP^N_{\alpha-\frac{1}{N}}
(t,z) \CP^N_{\beta-\frac{1}{N}}
(t,z)}{\left(\sum_{\alpha}
\CR_t^N(\alpha)\CP^N_{\alpha}(t,z)\right)^2} 
-\frac{\sum_{\alpha,\beta>0} (\alpha-\beta)^2
 \CP^N_{\alpha-\frac{1}{N}}
(t,z) \CP^N_{\beta-\frac{1}{N}}
(t,z)}{\left(\sum_{\alpha} \CP^N_{\alpha}(t,z)\right)^2} \right|\\
&\leq & C\left|\frac{\left(\sum_{\alpha,\beta>0}
(\alpha-\beta)^2 \CR_t^N(\alpha)\CR_t^N(\beta)
 \CP^N_{\alpha-\frac{1}{N}}
(t,z) \CP^N_{\beta-\frac{1}{N}}
(t,z)\right)\left(\sum_{\alpha} \CP^N_{\alpha}(t,z)\right)^2}{ \left( \sum_{\alpha,\beta>0}
(\alpha-\beta)^2
 \CP^N_{\alpha-\frac{1}{N}}
(t,z) \CP^N_{\beta-\frac{1}{N}}
(t,z)\right) \left(\sum_{\alpha}
\CR_t^N(\alpha)\CP^N_{\alpha}(t,z)\right)^2   }
 -1\right| .
\end{eqnarray*}

By localizing the summand and similar argument for the interior, we can show that for any $\epsilon' >0$, there exists a sufficiently large $N'$ such that for all $N\geq N'$,
$$\frac{1}{N}\left|\frac{\partial^2}{\partial r^2} \log
\ecal_N(t,z)-\frac{\partial^2}{\partial r^2} \log \Pi_N(t, z)\right|<\epsilon'.$$

\bigskip

\noindent {\bf Second order convergence in pure time derivatives}

\bigskip

We now consider time derivatives.  Let $G_i(x)=G_P(x)+f_i(x)$, for
$i=0, 1$. Let $f_t(x)=(1-t)f_0(x)+tf_1(x)$ and
$\frac{\partial}{\partial t}f_t(x)=v(x)=f_1(x)-f_0(x)$. Also $U_t(\rho) = u_t(z)$. By Legendre transform, $U_t(\rho)= x \rho - G_t(x)$ with $\rho = G_t'(x)$ and $x= U_t'(\rho)$. Calculate
$$\ddt{} U_t(\rho) = \dot x \rho - \ddt{} G_t(x)- G_t'(x) \dot x = -\ddt{} G_t(x) = - v(x)$$
and
$$\ddtt{} U_t(\rho) = -v'(x) \dot x =(v'(x))^2 U_t''(\rho)=(v'(x))^2\frac{\partial^2}{\partial \rho^2} u_t(z).$$

Straightforward calculation shows that

\begin{eqnarray*}
\frac{1}{N} \frac{\partial}{\partial t} \log
\ecal_N(t,z)+\ddt{} u_t(z)&=& \frac{1}{N} \frac{\sum_{\alpha}
\log\frac{Q_0^N(\alpha)}{Q_1^N(\alpha)} \CR_t^N(\alpha)
\CP^N_{\alpha}(t,z)}{\sum_{\alpha}\CR_t^N(\alpha)\CP^N_{\alpha}(t,z)} \\
&=& \frac{1}{N} \frac{\sum_{\alpha} \left(-N
v(\alpha)+\log\frac{q_0^N(\alpha)}{q_1^N(\alpha)}\right)
\CR_t^N(\alpha)\CP^N_{\alpha}(t,z)}{\sum_{\alpha}\CR_t^N(\alpha)
\CP^N_{\alpha}(t,z)}
\end{eqnarray*}
and
\begin{eqnarray*}
&&\frac{1}{N}\frac{\partial^2}{\partial t^2} \log \ecal_N(t,z) +  \ddtt{} u_t(z)\\
&=& N \frac{\sum_{\alpha,\beta} \left((v(\alpha)-v(\beta))
+\frac{1}{N} \log\frac{q_0^N(\beta)q_1^N(\alpha)}
{q_0^N(\alpha)q_1^N(\beta)}\right)^2 \CR_t^N(\alpha)\CR_t^N(\beta)
\CP^N_{\alpha}(t,z)\CP_t^N(\beta,z)}{ 2 \left(\sum_{\alpha}\CR_t^N(\alpha)
\CP^N_{\alpha}(t,z)\right)^2}\\
&=& I_1^N(t,z)+I_2^N(t,z)+I_3^N(t,z).
\end{eqnarray*}
\begin{eqnarray*}
&&I_1^N(t,z)\\
&=& N \frac{\sum_{\alpha,\beta} (v(\alpha)-v(\beta))^2 \CR_t^N(\alpha)\CR_t^N(\beta)
\CP^N_{\alpha}(t,z)\CP_t^N(\beta,z)}{2 \left(\sum_{\alpha}\CR_t^N(\alpha)
\CP^N_{\alpha}(t,z)\right)^2}\\
&\sim& N \frac{\sum_{\alpha,\beta \in B_{x}(\delta) } (\alpha-\beta)^2 \left( \frac{v(\alpha)-v(\beta)}{\alpha-\beta}\right)^2 \CR_t^N(\alpha)\CR_t^N(\beta)
\CP^N_{\alpha}(t,z)\CP_t^N(\beta,z)}{ 2 \left(\sum_{\alpha}\CR_t^N(\alpha)
\CP^N_{\alpha}(t,z)\right)^2}\\
&\sim& N(v'(x))^2 \frac{\sum_{\alpha,\beta \in B_{x}(\delta) } (\alpha-\beta)^2\CR_t^N(\alpha)\CR_t^N(\beta)
\CP^N_{\alpha}(t,z)\CP_t^N(\beta,z)}{2\left(\sum_{\alpha}\CR_t^N(\alpha)
\CP^N_{\alpha}(t,z)\right)^2}\\
&&+N \frac{\sum_{\alpha,\beta \in B_{x}(\delta) } (\alpha-\beta)^2  \left( \left( \frac{v(\alpha)-v(\beta)}{\alpha-\beta}\right)^2 -(v'(x))^2\right) \CR_t^N(\alpha)\CR_t^N(\beta)
\CP^N_{\alpha}(t,z)\CP_t^N(\beta,z)}{ 2 \left(\sum_{\alpha}\CR_t^N(\alpha)
\CP^N_{\alpha}(t,z)\right)^2}\\
&\sim& \left((v'(x))^2 + \sup_{\alpha\neq \beta \in B_{x}(\delta)}  \left( \left( \frac{v(\alpha)-v(\beta)}{\alpha-\beta}\right)^2 -(v'(x))^2\right) \right) \left(\frac{1}{N}\frac{\partial^2}{\partial \rho^2} \log \ecal_N(t,z) + \frac{\partial^2}{\partial \rho^2} u_t(z)\right)\\
&\sim& (v'(x))^2 \frac{\partial^2}{\partial \rho^2} u_t(z)=\ddtt{} u_t(z)
\end{eqnarray*}
as $\delta\rightarrow 0$.

Therefore $\lim_{N\rightarrow \infty} I_1^N(t,z) = \ddtt{} u_t(z)$.

\begin{eqnarray*}
I_2^N(t,z)&=& \frac{\sum_{\alpha,\beta} (v(\alpha)-v(\beta))
 \log\frac{q_0^N(\beta)q_1^N(\alpha)}
{q_0^N(\alpha)q_1^N(\beta)} \CR_t^N(\alpha)\CR_t^N(\beta)
\CP^N_{\alpha}(t,z)\CP_t^N(\beta,z)}{\left(\sum_{\alpha}\CR_t^N(\alpha)
\CP^N_{\alpha}(t,z)\right)^2}\\
& \sim & \frac{\sum_{\alpha,\beta\in B_x(\delta)} (v(\alpha)-v(\beta))
 \log\frac{q_0^N(\beta)q_1^N(\alpha)}
{q_0^N(\alpha)q_1^N(\beta)} \CR_t^N(\alpha)\CR_t^N(\beta)
\CP^N_{\alpha}(t,z)\CP_t^N(\beta,z)}{\left(\sum_{\alpha}\CR_t^N(\alpha)
\CP^N_{\alpha}(t,z)\right)^2}\\
&\sim& 0
\end{eqnarray*}
as $\delta\rightarrow 0$.

Therefore $\lim_{N\rightarrow \infty} I_2^N(t,z) = 0$.

\begin{eqnarray*}
I_3^N(t,z)&=& \frac{1}{N}\frac{\sum_{\alpha,\beta} \left(
 \log\frac{q_0^N(\beta)q_1^N(\alpha)}
{q_0^N(\alpha)q_1^N(\beta)} \right)^2 \CR_t^N(\alpha)\CR_t^N(\beta)
\CP^N_{\alpha}(t,z)\CP_t^N(\beta,z)}{\left(\sum_{\alpha}\CR_t^N(\alpha)
\CP^N_{\alpha}(t,z)\right)^2}\\
&\sim& 0.
\end{eqnarray*}

Therefore $\lim_{N\rightarrow \infty} I_3^N(t,z) = 0$.

We conclude from the above calculation that
$$\lim_{ N \rightarrow \infty}  \ddtt{} \log \ecal_N(t,z) = 0.$$

By a similar argument, which we leave to the reader,  the mixed
space-time derivatives of $\log \mathcal{E}_N(t, z)$  also
uniformly converges to $0$. Therefore  $\frac{1}{N} \log \mathcal{E}_N(t,
z)$ has bounded second derivatives and $\frac{1}{N} \log
\mathcal{E}_N(t, z)$ converges in $C^{2}$ to $0$.
\qed
\end{proof}

We now conclude the proof of the main result:

\bigskip

\noindent {\bf Proof of Theorem \ref{main}}

\bigskip

\noindent Notice that $\varphi_N(t, \cdot)-\varphi_t(\cdot)=\frac{1}{N} \log \ecal_N (t, \cdot)$.
Therefore Theorem \ref{main} is proved.  \qed

\bigskip

\subsection{Final remarks and questions}

We conclude with some questions:

\begin{itemize}

\item Limits on the degree of convergence of $\varphi_N(t,z) \to
\varphi(t,z) $ are related to the distribution of complex zeros of
the holomorphic extension of $\ecal_t^N$. In the toric case, in
the coordinates $x = e^{\rho}$, $\ecal_t^N$ is a positive real
polynomial of a real variable. As observed by  Lee-Yang in the
context of partition functions of statistical mechanical models,
the degree of convergence of $\frac{1}{N} \log \ecal_t^N$ to its
limit is related to the limit distribution  of the complex zeros
of $\ecal_t^N$ along the real domain. It would be interesting to
study the complex zeros in the case of toric varieties. The
approximation (\ref{BCALK}) is in some sense a `best polynomial
approximation' to a general \kahler potential, and it is possible
that there exist analogues of results relating convergence of such
approximations to distribution of zeros.

\item The formula for $\ecal_t^N(z)$ in Lemma \ref{RCAL} exhibits
this function as the value on the diagonal of a Toeplitz type
operator with multiplier $\rcal_t^N(\alpha)$. More precisely, it
is the Berezin lower symbol of the Toeplitz type operator. For
background we refer to \cite{STZ2}.  The question whether it is a
Toeplitz operator in any standard sense is essentially the same
question as to the existence of asymptotics of $\ecal_t^N(z)$ and
joint asymptotics of $\rcal_t^N(\alpha)$. It would be very helpful
if there exists a more `abstract' approach to this Toeplitz
operator by constructing its Toeplitz symbol instead of its
Berezin symbol. The leading order Toeplitz symbol is calculated in
Corollary \ref{maincor1}.

\end{itemize}


\end{document}